\documentclass[letterpaper, 10 pt, conference]{ieeeconf}
\usepackage{amssymb}
\usepackage{color}
\usepackage{bbm}

\usepackage{graphicx}
\usepackage{mathtools}
\usepackage{multirow}
\usepackage{lineno,hyperref}
\usepackage{url}
\usepackage{cleveref}
\usepackage{subcaption}
\usepackage{bm}
\usepackage{tikz}
\usetikzlibrary{arrows.meta}
\usepackage{algorithm}
\usepackage{algorithmic}
\usepackage{caption}
\usepackage{mathtools}
\mathtoolsset{showonlyrefs,showmanualtags}

\DeclareMathAlphabet{\pazocal}{OMS}{zplm}{m}{n}

\newcommand{\ve}{\boldsymbol}

\renewcommand{\epsilon}{\varepsilon}

\usepackage{pgfplots}
\DeclareUnicodeCharacter{2212}{−}
\usepgfplotslibrary{groupplots,dateplot}
\usetikzlibrary{patterns,shapes.arrows}

\modulolinenumbers[5]

\usepackage{multirow}
\usepackage{color}
\usepackage{paralist}
\usepackage[space]{cite}

\pgfplotsset{compat=1.18}

\IEEEoverridecommandlockouts                        

\title{\LARGE \bf
Sample-Efficient Counterfactual Tuning for Compressor Pressure Control
}

\author{Margarita A. Guerrero$^{1}$, Rodrigo A. González$^{2}$,
and Cristian R. Rojas$^{1}$%
\thanks{This work has been partially supported by the Swedish Research Council under contract number 2023-05170, and by the Wallenberg AI, Autonomous Systems and Software Program (WASP) funded by the Knut and Alice Wallenberg Foundation.}%
\thanks{$^{1}$Margarita A. Guerrero and Cristian R. Rojas are with the
Department of Decision and Control Systems, KTH Royal Institute of
Technology, 100 44 Stockholm, Sweden
(e-mails: mags3@kth.se, crro@kth.se).}%
\thanks{$^{2}$Rodrigo A. González is with the Department of Mechanical
Engineering, Eindhoven University of Technology,
Eindhoven, The Netherlands
(e-mail: r.a.gonzalez@tue.nl).}%
\thanks{Code available at \url{https://github.com/mags-ono/ACE}}
}

\begin{document}

\maketitle

\begin{abstract}
In controlled industrial environments, ensuring safety and performance during controller tuning is a challenging and critical task. In particular, control loops in compressor–plenum–throttle systems cannot tolerate costly interruptions, and aggressive excitation may lead to unsafe operating regimes. Given the wide availability of historical data under safe operation, this paper introduces a counterfactual explainability approach for sample-efficient retuning of compressor control loops. The proposed data-driven algorithm determines, without an explicit plant model or previous control law, the smallest controller adjustment required to achieve predefined performance specifications while maintaining closed-loop stability. The effectiveness of the method is demonstrated through an extensive Monte Carlo simulation study, where the proposed approach computes low magnitude and interpretable adjustments in the controller space that lead to improved settling time and reduced overshoot.
\end{abstract}



\section{Introduction}

In industrial facilities, not all units carry the same operational weight; rather, a select group of critical assets governs whether the plant remains on-spec and online, shaping availability, quality, and safety downstream. Among these, compressors stand out for their cross-sector impact: in natural-gas transmission, compressor stations maintain pressure and flow over hundreds of kilometers of pipelines and operate continuously as bottleneck assets \cite{kurz2012gas}. In natural gas liquefaction, large refrigeration trains are dominated by compressor power and availability, 
while in air-separation units compression underpins cryogenic distillation and core safety practices~\cite{smith2001}. Finally, from grid-scale gas-turbine power plants to wide-body aero engines~\cite{brandt1994}, axial compressors are core hardware that must operate within a safe envelope. Across these process industries, that envelope is bounded by rotating stall, a localized aerodynamic loss, and surge, a global instability with pressure and flow oscillations, which must be avoided to ensure stable operation~\cite{greitzer1976surge}.

Control loops around compressors are critical to maintain stable flow and adequate pressurization; however, tuning these loops is notoriously hard in practice~\cite{jager95}. Opening the loop to identify the plant or to run relay tests is usually infeasible: the unit is a bottleneck, interruptions are costly, and aggressive excitation risks flow reversal. Even during maintenance, the load seen by the unit 
is rarely representative of normal operation. 

A natural prior for closed-loop tuning is the plant’s Distributed Control System (DCS) or Programmable Logic Controller (PLC) historian~\cite{crompton2021}, a time-series database that typically archives long-horizon records of the process variable (PV), setpoint (SP), controller output (OP), and controller mode, supporting envelope characterization, representative-window selection, and time-domain benchmarking near the operating point. However, two practical obstacles arise. First, the controller numerical parameter values (e.g., PID gains or lead--lag time constants) are often not logged, thus parameter identification is required. Second, historical segments must be representative of a time-invariant plant; otherwise, observed changes reflect plant drift 
rather than controller action, biasing identification and tuning. These realities motivate historian-anchored priors combined with minimal, safe, closed-loop tests near the current operating point.

Within this landscape, the axial-compressor pressure loop, regulated via Inlet Guide Vanes (IGVs)---adjustable inlet blades that modulate mass flow and incidence---offers an ideal testbed. Their unsteady behavior near healthy operating points admits compact physics–based reductions such as the Moore-Greitzer compressor-plenum-throttle model \cite{moore1985theory}, while under correct operation and standard maintenance the hardware envelope remains remarkably stable over multi-year horizons~\cite[p.~29]{beagle21}, making DCS historian records comparable across long spans.  However, safe retuning is tightly constrained by surge avoidance and unit availability: only short, closed-loop tests are acceptable, and aggressive excitation is off–limits.

Against this backdrop, there is a substantial body of work on data-driven controller design, including closed-loop identification with iterative redesign~\cite{anderson1991adaptive,van1995identification,albertos2012iterative}, model-free iterative tuning~\cite{hjalmarsson1998iterative,campi2002virtual}, and relay autotuning~\cite{aastrom2006advanced}. While powerful in general settings, these approaches are primarily aimed at achieving desired control performance via dedicated excitation experiments or model-based redesign procedures, and they provide only limited guarantees that the controller update will be both small and operator-actionable (i.e., restricted to permitted parameters within tight bounds). Moreover, they face practical limitations in safety-critical compressor loops, as they typically require non-trivial excitation or operating-point perturbations that are unacceptable near surge margins.

We address the controller retuning problem by posing it as a Counterfactual Explainability (CFE) problem~\cite{wachter_17}. A counterfactual explanation constitutes the smallest, actionable change in the inputs that flips a model's outcome.  We cast controller retuning exactly in that form—find the minimal, operator-actionable change in the controller parameters that makes the pressure-control loop satisfy time-domain specifications—while respecting tight excitation constraints and a trust region around the baseline tuning. In this sense, CFEs provide an action-oriented lens within explainable AI (XAI)~\cite{Molnar_22}, focusing not just on why a decision is made but on how to change it with the least intervention. Within XAI for control, recent work has adapted local attribution tools such as LIME~\cite{Ribeiro16} and SHAP~\cite{Lundberg17} to interpret controller behavior and obtain insights into optimization choices for closed-loop performance~\cite{allamaa2025,porcari2025}, while property-preserving XAI provides explanations that retain control-relevant guarantees~\cite{Riva24}. In contrast to these descriptive approaches, we take an action-oriented route: inspired by Adaptive Sampling for Counterfactual Explanations (ACE)~\cite{Guerrero-25}, we compute the minimal change to controller parameters that achieves the required closed-loop specifications with as few experiments as possible, restricted to a trust region around the current operating point.
\medskip

The key idea is to cast “meeting the time-domain specs” as a binary classification task for the pressure controller that regulates the pressure-rise coefficient $\Psi$—i.e., the atmosphere-to-plenum pressure difference—via IGV actuation, while maintaining a safe surge margin. Starting from a controller that does not meet the closed-loop specifications, we then search the controller-parameter space for the nearest label-flipping modification that turns “fail” into “pass,” staying close to the current setting and within a safe neighborhood. Historian data provide a prior over feasible behavior; a Gaussian-Process (GP) classifier supplies probabilistic pass scores; and a Bayesian-optimization loop proposes a few targeted, safe experiments near the data manifold to validate and refine the retuning.

To the best of our knowledge, this is the first formulation of adaptive sampling for counterfactual explanation tailored for control tuning. Our main contributions are:
\begin{itemize}
\item A historian-anchored, CFE formulation of safe, minimum change controller retuning for compressor loops.
\item An ACE-style Bayesian procedure that requires only a handful of closed-loop step tests, guided by a GP classifier and trust-region acquisitions.
\item We demonstrate the effectiveness of the proposed controller tuning algorithm on a Compressor–Plenum– Throttle system, highlighting how the critical pressure control loop meets specifications by searching within a safe region, avoiding surge.
\end{itemize}







The remainder of the paper is structured as follows. Section \ref{sec:processdescription} describes the compressor system and formulates the controller retuning problem. Section \ref{sec:cfe} contains the main contributions of this work, namely, a counterfactual explainability approach to controller retuning in compressor loops. Simulation results showcasing the effectiveness of the proposed approach are presented in Section \ref{sec:results}, while conclusions are drawn in Section \ref{sec:conclusions}.

\section{Process description and problem formulation}
\label{sec:processdescription}
In this section, we present the canonical air system, an axial compressor discharging into a downstream plenum, and the closed-loop architecture used for controller retuning.

\subsection{Physical setting}
In industrial process and turbomachinery applications, compressors raise total pressure to deliver the target mass flow against downstream demand. In axial compressors, rotor–stator rows add work and guide the flow, while IGVs set inlet pre-swirl to adjust incidence and mass flow at a fixed shaft speed, shifting the operating point along the compressor characteristic $\psi_c$ (the quasi-steady map of pressure-rise coefficient versus flow coefficient at fixed speed and IGV setting).

\begin{figure}[h!]
\centering
\includegraphics[trim=0mm 0mm 0mm 0mm,clip, width=1\columnwidth]{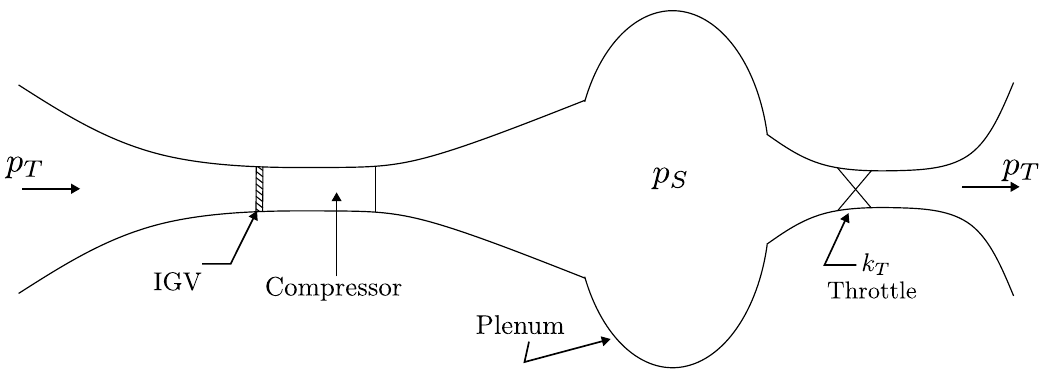}
\caption{Schematic of the compressor system geometry.}
    \label{fig:compressor_system} 
\end{figure}

Within this framework, we consider an axial compressor installed in a duct and discharging into a downstream plenum (Fig.~\ref{fig:compressor_system}). The plenum volume is large relative to the compressor duct; hence, flow speeds and fluid accelerations in the plenum are negligible, and the plenum pressure $p_S$ can be treated as spatially uniform, albeit time-varying. We denote by $p_T$ the upstream (ambient/atmospheric) pressure. A throttle at the plenum exit represents the downstream consumer (e.g., a gas turbine) and imposes a linear pressure–flow slope $k_T$. The compressor IGVs are actuated by a servo valve; by modulating the IGV position $u$, the operator sets the pressure rise that sustains the downstream unit. This lumped description is standard in the Moore–Greitzer model~\cite{moore1985theory} and captures the key interactions of the compressor–plenum–throttle system.

A key operational constraint is \emph{surge}: if plenum pressure rises excessively---i.e., the atmosphere-to-plenum pressure rise $\Psi$ increases beyond the compressor’s capability---the annulus-averaged flow may reverse, driving pressure oscillations that threaten both blades and downstream hardware. Our study focuses on a healthy operating region---far from surge---while acknowledging that pressure excursions can move the system toward this boundary. In industrial practice, the hardware and its operating envelope remain remarkably stable over multi-year periods, thus archived process data provide a suitable prior information for identification and controller retuning.

Under the Moore–Greitzer assumptions---axisymmetric, annulus-averaged flow; a quasi-steady compressor characteristic $\psi_c(\Phi)$ at fixed shaft speed and IGV setting; and the physical setting of Fig.~\ref{fig:compressor_system}---the compressor–plenum system reduces to a second-order system for the annulus-averaged axial-flow coefficient $\Phi$:
\begin{equation}\label{eq:MG46}
\resizebox{0.48\textwidth}{!}{$
\frac{d^{2}\Phi}{d\xi^{2}}
+ \frac{1}{\ell_{c}}\!\left(\frac{1}{4B^{2}k_{T}} - \frac{d\psi_{c}}{d\Phi}\right)\!\frac{d\Phi}{d\xi}
+ \frac{1}{4\ell_{c}^{2}B^{2}}
\left[(\Phi - \bar{\Phi}) - \frac{\Psi - \bar{\Psi}}{k_{T}}\right] = 0,
$}
\end{equation}
where $\psi_c$ is the compressor characteristic, $\Psi$ the pressure rise coefficient proportional to $p_S - p_T$ (plenum relative to atmosphere), $\xi=Ut/R$ the nondimensional time with $U$ the rotor peripheral speed at the mean radius $R$, and $u$ the IGV position (per unit). The constants are $\ell_c=L_c/R$ (duct-length parameter), $B$ (Greitzer parameter), and $k_T$ (throttle slope). All variables are nondimensional, and overbars (e.g., $\bar{\Phi}$, $\bar{\Psi}$) denote operating-point values. Equation \eqref{eq:MG46} is the well-known Moore–Greitzer pure-surge reduction, highlighting how $B$ influences both the damping and stiffness of the lumped dynamics. 

\subsection{IGV control}

The pressure-rise from atmosphere to plenum must be controlled in a feedback loop. To this end, starting from the Moore-Greitzer dynamics in \eqref{eq:MG46}, a small-signal model is derived around a healthy operating point under IGV actuation. Linearizing around $(\bar{\Phi},\bar{u},\bar{\Psi})$ far from surge and defining
$\tilde{\Phi}=\Phi-\bar{\Phi}$, $\tilde{u}=u-\bar{u}$, and $\tilde{\Psi}=\Psi-\bar{\Psi}$, the compressor map is locally
\begin{equation}
\psi_c(\Phi,u)\;\approx\;\psi_c(\bar{\Phi},\bar{u})\;+\;a\,\tilde{\Phi}\;+\;\alpha\,\tilde{u},
\label{eq:small-signal-map}
\end{equation}
with $a=\frac{\partial\psi_c}{\partial\Phi}\big|_{(\bar{\Phi},\bar{u})},\;
\alpha=\frac{\partial\psi_c}{\partial u}\big|_{(\bar{\Phi},\bar{u})}$. 
The linearized momentum equation in \cite[Chap.~4, Eq.~(4.2)]{moore1985theory} yields the small-signal relation
\begin{equation}
\tilde{\Psi}
= a\,\tilde{\Phi}
+ \alpha\,\tilde{u}
- \ell_c\,\frac{d\tilde{\Phi}}{d\xi}.
\label{eq:psi-phi-u}
\end{equation}
Combining \eqref{eq:MG46} with \eqref{eq:small-signal-map}–\eqref{eq:psi-phi-u} and taking Laplace transforms gives a second-order transfer function $G_1(s)$ from IGV to pressure rise,
\[
G_1(s)= \frac{c_1\,s + c_0}{s^{2} + 2\zeta\omega_n s + \omega_n^{2}},
\]
with $s$ being the Laplace transform variable, and $\omega_n$, $\zeta$, $c_0$, and $c_1$ set by $(\ell_c,B,k_T,a,\alpha)$. 

\begin{figure}[h!]
\centering
\includegraphics[trim=0mm 0mm 0mm 0mm,clip, width=1\columnwidth]{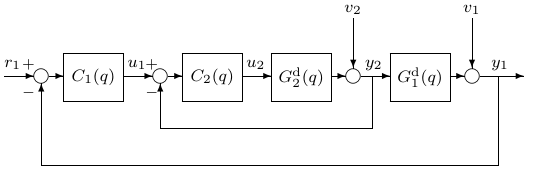}
	\caption{Cascaded system description for IGV control, with $C_1(q)$, $C_2(q)$, $G_1^{\text{d}}(q)$ and $G_2^{\text{d}}(q)$ being the outer controller, inner controller, compressor-plenum system, and IGV actuator, respectively.}
	\label{fig_cascade}
\end{figure}

On the other hand, the IGV actuator is modeled as a first-order lag with static gain $k_v$ and time constant $T_v$, a standard approximation for electro–hydraulic servo valves:
\[
G_2(s) \;=\; \frac{k_v}{1 + T_v s}.
\]
The dynamics of the complete compressor-\allowbreak plenum-\allowbreak throttle system are controlled by a discrete-time cascaded control structure operating at the sampling frequency set by the DCS (see Fig.~\ref{fig_cascade}). Under this architecture the continuous-time transfer functions $G_1(s)$ and $G_2(s)$ are replaced by their discrete-time equivalents \cite{aastrom1984computer} $G_1^{\textnormal{d}}(q)$ and $G_2^{\textnormal{d}}(q)$, where $q$ denotes the forward shift operator. The blocks $C_1(q)$ and $C_2(q)$ are designable linear and time-invariant controllers, and $v_1$ with $v_2$ are stationary stochastic processes that model output disturbances and sensor noise. In addition, $y_1$ is the pressure-related controlled variable (e.g., normalized pressure rise), and $r_1$ is the reference pressure rise. In practice, $r_1,u_1,u_2,y_1$ and $y_2$ are typically available from the DCS historian, whereas $v_1,v_2,G_1^{\textnormal{d}}(q)$ and $G_2^{\textnormal{d}}(q)$ are not.

\subsection{Problem formulation}\label{sec:problem_formulation}

This paper addresses the problem of controller retuning for critical feedback loops in axial-compressor (compressor-\allowbreak plenum-\allowbreak throttle) systems. The goal is to derive an iterative procedure that produces the smallest possible change in the baseline controller parameters while maintaining safe operation and ensuring control performance. 

The controller is parameterized by a vector $\boldsymbol{\theta}\in\Theta$ (for example, PID gains), with $\boldsymbol{\theta}_0$ denoting the baseline control law. We assume the availability of $N_{\textnormal{hist}}$ batches of historical data $\{r_1^{(i)}[k],u_{1}^{(i)}[k],u_{2}^{(i)}[k]$, $y_{1}^{(i)}[k],y_{2}^{(i)}[k]\}_{k=1}^{N_i}$, $i=1,\dots, N_{\textnormal{hist}}$, collected under distinct control laws. Since controller gains are not logged, they must be estimated from the data. The models $G_1^{\textnormal{d}}(q)$ and $G_2^{\textnormal{d}}(q)$ represent operation in a healthy region away from surge, and, because the hardware and operating envelope are remarkably stable over multi-year horizons, the plants are assumed to be time-invariant across the batches.

We seek a controller, parameterized by $\boldsymbol{\theta}$, that remains as close as possible to $\boldsymbol{\theta}_0$ under a prescribed metric $d$ on $\Theta$, while satisfying time-domain performance metrics such as fast response, stability, bounded overshoot, and small steady-state error. This leads to the optimization problem
\begin{equation}
\displaystyle \min_{{\boldsymbol{\theta}} \in \Theta} 
d({\boldsymbol{\theta}}, \boldsymbol{\theta}_0) 
\quad \text{s.t.} \quad 
\textnormal{performance specifications met}.
\label{eq:optim_problem}
\end{equation}

Because experiments on real compressor systems are costly and safety critical, only closed-loop data is collected, and the number of experiments must be kept minimal. 

\section{CFE approach to sample-efficient control tuning}
\label{sec:cfe}
The objective in \eqref{eq:optim_problem} is naturally cast as a counterfactual-explanation task on the controller space: starting from the factual $\boldsymbol{\theta}_0$ (label ``fail''), and the available historical data, find the nearest counterfactual---i.e., a solution to \eqref{eq:optim_problem}---that flips the label to ``pass''. This formulation is well aligned with closed-loop tuning, since evaluating a candidate $\ve{\theta}$ yields a
binary outcome $h(\ve{\theta})\in\{0,1\}$ that provides an intuitive indicator of specification compliance for retuning. In the sequel we operationalize this idea with an ACE-style, Bayesian optimization scheme that (i) leverages the DCS historian to build a prior on feasible controllers around historically feasible behavior and (ii) searches within a trust region near the current operating point, thereby proposing safe minimal-change controller updates without opening the loop.

\subsection{Controller identification and performance evaluation}
\label{sec:controller}
A key step in the proposed approach is the collection of a labeled dataset from historical closed-loop data. Specifically, the controller parameters $\boldsymbol{\nu}_1^{(i)}$ and $\boldsymbol{\nu}_2^{(i)}$ for $C_1(q)$ and $C_2(q)$ must be estimated, given the $i$th batch of cascaded system data in Fig. \ref{fig_cascade}, $(\boldsymbol{r}_1^{(i)}, \boldsymbol{u}_1^{(i)}, \boldsymbol{u}_2^{(i)}, \boldsymbol{y}_1^{(i)}, \boldsymbol{y}_2^{(i)})$, where $\boldsymbol{r}_1^{(i)}=\{r_1[k]\}_{k=1}^{N_i}$ and analogously for the other signals. Methods for this identification step include instrumental variables, prediction error methods, or errors-in-variable methods \cite{Ljung1999SystemUser}. Since both controller input and output signals are assumed to be measured, a straightforward approach that ensures reliable estimates is to solve the corresponding nonlinear least-squares problems:
\begin{align}
    \boldsymbol{\nu}_1^{(i)} &= \underset{\boldsymbol{\nu}_1}{\arg\min} \hspace{0.1cm} \sum_{k=1}^{N_i} \left(u_{1}^{(i)}[k]-C_1(q,\boldsymbol{\nu}_1) (r_{1}^{(i)}[k]-y_{1}^{(i)}[k])\right)^2, \notag \\
    \boldsymbol{\nu}_2^{(i)} &= \underset{\boldsymbol{\nu}_2}{\arg\min} \hspace{0.1cm} \sum_{k=1}^{N_i} \left(u_{2}^{(i)}[k]-C_2(q,\boldsymbol{\nu}_2) (u_{1}^{(i)}[k]-y_{2}^{(i)}[k])\right)^2. \notag
\end{align}
The quality of the identification step depends on the richness of the data and the sample size. The reference signal must be persistently exciting of sufficiently high order, with step tests designed to provide rich data according to standard input design rules \cite{Soderstrom2001SystemIdentification}. In addition, in the presence of noise affecting the measurements, the asymptotic covariance of the estimated parameters decreases proportionally to $1/N_i$ \cite{Ljung1999SystemUser}.

Let $\boldsymbol{\theta}^{(i)} = [\boldsymbol{\nu}_{1}^{(i)\top}, \boldsymbol{\nu}_{2}^{(i)\top}]^\top \in \mathbb{R}^{p}$ denote the controller parameters pertaining to the $i$th batch of data. Given a closed-loop step experiment, we evaluate a vector of time-domain specifications
\[
\pi(\boldsymbol{\theta}^{(i)}) \;=\; \bigl(e_{\mathrm{ss}}(\boldsymbol{\theta}^{(i)}),\, t_{\mathrm{rise}}(\boldsymbol{\theta}^{(i)}),\, t_{\mathrm{settle}}(\boldsymbol{\theta}^{(i)}),\, \mathrm{OS}(\boldsymbol{\theta}^{(i)})\bigr),
\]
where $e_{\mathrm{ss}}$ is the steady-state error, $t_{\mathrm{rise}}$ is the rise time, $t_{\mathrm{settle}}$ is the settling time for a fixed amplitude margin, and OS is the overshoot percent. Given $\pi(\boldsymbol{\theta}^{(i)})$ for each $i=1,\dots,N_{\textnormal{hist}}$, we define the binary labels
\[
h(\boldsymbol{\theta}^{(i)}) \;=\; \mathbb{I}\!\{ \pi(\boldsymbol{\theta}^{(i)}) \preceq \tau \},
\quad
 i = 1,\dots, N_{\textnormal{hist}},
\]
where $\tau=\bigl(\tau_{e},\,\tau_{\mathrm{rise}},\,\tau_{\mathrm{settle}},\,\tau_{\mathrm{OS}}\bigr)$ describes the performance specifications for each metric. Note that $h(\boldsymbol{\theta}^{(i)})$ equals $1$ if and only if all specifications are satisfied. The oracle $h(\cdot)$ is a black box (obtained via closed-loop simulation/operation) and costly to query; safety further requires that queries remain near historically observed, feasible behavior, since available closed-loop data (collected under healthy operation, away from surge) define
the admissible operating envelope, i.e., the range of operating conditions and controller settings observed to be safe. All historical data is thus condensed in the labeled dataset $\mathcal{D}_{\text{hist}} = \{ (\boldsymbol{\theta}^{(i)}, \ell_i) \}_{i=1}^{N_{\text{hist}}}$, where $\ell_i = h(\boldsymbol{\theta}^{(i)}) \in \{0,1\}.$



\subsection{Counterfactual explainability step}
Having identified controllers from historical data and assembled a labeled dataset $\mathcal{D}_{\text{hist}}$, 
we now use this information to construct the adaptive sampling procedure for counterfactual explanations.
Given a current controller $\ve \theta_0$ that fails the time-domain specifications (i.e., $h(\ve\theta_0)=0$), the goal is to find a new controller \(\ve{\theta}\) that satisfies the specifications while remaining as close as possible to $\ve \theta_0$.

Following the Gaussian process classification method \cite{Rasmussen-Williams-06}, we rewrite the pass/fail rule from Section~\ref{sec:controller} as
\[h(\ve{\theta})=\mathbb{H}(f(\ve{\theta})-0.5),\]
where $\mathbb{H}\colon \mathbb{R} \to \{0, 1\}$ is the Heaviside step function (i.e., $\mathbb{H}(a)=1$ if $a\ge 0$ and $0$ otherwise), and $f \colon \Theta \to \mathbb{R}$ is a smooth latent function. We refer to $f$ as the \textit{black-box target function}, where the latent score $f(\ve\theta)$ summarizes how the closed-loop step response under controller $\ve\theta$ compares against the time-domain specifications: values above $0.5$ indicate “pass” and below $0.5$ indicate “fail”.

At a high level, using historical data we obtain prior knowledge over feasible behavior, then adaptively propose new controllers near the data manifold and the current setting $\boldsymbol{\theta}_0$. Each proposal triggers a short closed-loop step test; we record the pass/fail label and restrict subsequent queries to a trust region to keep trials safe and minimal. Thus, the historian anchors the prior, and ACE’s acquisition steers a few targeted experiments toward the nearest label-flipping retuning that meets the specifications.

In what follows, we instantiate ACE for compressor–IGV tuning, following the guidelines presented in \cite{Guerrero-25}.



\subsection*{Cost function}
The ACE method 
interprets the decision boundary as the set $\{\,\boldsymbol{\theta}\in \Theta\colon f(\boldsymbol{\theta})=0.5\,\}$, i.e., where pass and fail are equally likely. Hence, the counterfactual problem in \eqref{eq:optim_problem} can be re-written as
\begin{equation} 
\underset{\ve \theta}{\text{min}} \ d(\ve \theta, \ve \theta_0) \quad \text{s.t.} \ f(\ve \theta) = 0.5. 
\label{eq:penal}
\end{equation}
Following \cite{Guerrero-25}, we adopt a penalty formulation to obtain an unconstrained program:
\begin{equation}\label{eq:EI_penalty}
\underset{\boldsymbol{\theta}}{\text{min}}\;\;
J(\boldsymbol{\theta}) \;=\;
d(\boldsymbol{\theta},\boldsymbol{\theta}_0)
\;+\; \lambda_k |f(\ve \theta)-0.5|,
\end{equation}
where $\lambda_k>0$ can be increased across iterations $k \in \mathbb{N}$ (denoting the ACE iteration index, cf. Algorithm.~\ref{alg:ace_dcs} in Section~\ref{sec:algo}) to tighten the boundary constraint. The terms $d(\boldsymbol{\theta},\boldsymbol{\theta}_0)$ and $|f(\boldsymbol{\theta})-0.5|$ jointly enforce \emph{proximity}: the first term keeps the candidate controller close to the current one in the parameter space, while the second one drives queries toward the decision boundary. Enforcing the boundary yields label changes with minimal parameter movement, consistent with the counterfactual goal.

To encode additional structural preferences, we add three terms—\emph{actionability} (search only over operator–changeable parameters), \emph{sparsity} (alter as few parameters as possible), and \emph{plausibility} (prefer controllers supported by historical data):
\begin{equation}\label{eq:asce_cost}
\resizebox{0.43\textwidth}{!}{$
\underbrace{\underset{\boldsymbol{\theta}\in\mathcal{A}\subset\Theta}{\arg\min}\;}_{\text{actionability}}
\underbrace{d(\boldsymbol{\theta}, \boldsymbol{\theta}_0) + \lambda_k\,\bigl|\,f(\boldsymbol{\theta})-0.5\,\bigr|}_{\text{proximity}}
\;+\;
\underbrace{\beta\, g(\boldsymbol{\theta}-\boldsymbol{\theta}_0)}_{\text{sparsity}}
\;+\;
\underbrace{l(\boldsymbol{\theta};\,\mathcal{D})}_{\text{plausibility}}
$}
\end{equation}

where:
(i) $\mathcal{A}\subset\Theta$ enforces \emph{actionability} by fixing non-actionable parameters (e.g., keeping the derivative gain in a PID controller unchanged when only retuning the remaining parameters);
(ii) $g(\boldsymbol{\theta}-\boldsymbol{\theta}_0)$ promotes \emph{sparsity} of changes (in this work, an $\ell_1$ norm);
(iii) $l(\boldsymbol{\theta};\mathcal{D})$ is a \emph{hard plausibility filter} implemented via the Local Outlier Factor
(LOF)~\cite{breunig_00}, assigning zero cost to inliers and discarding outliers far from the empirical data manifold, computed over $\mathcal{D}=\mathcal{D}_{\mathrm{hist}}\cup\mathcal{D}_{\mathrm{online}}$ (historical DCS data plus newly collected
closed-loop samples). As in \eqref{eq:EI_penalty}, the penalty method increases $\lambda_k$ across iterations; $\beta>0$ trades off sparsity against proximity and plausibility, and $d(\cdot,\cdot)$ is the Euclidean distance.

To this end, we define the \emph{trust region} as the set of admissible candidates explored by ACE around
$\boldsymbol{\theta}_0$, constrained to remain within historically observed feasible behavior (initial and newly
collected healthy data) and reinforced by the plausibility filter $l(\boldsymbol{\theta};\mathcal{D})$ in
\eqref{eq:asce_cost}.

\subsection*{Surrogate model}
Recall that $f$ is unknown and we only observe a binary label
$\ell = h(\boldsymbol{\theta})\in\{0,1\}$. Therefore, to solve \eqref{eq:asce_cost} we approximate the true underlying function $f$ with a surrogate $\hat f$. In our setting, $\hat f$ is a Gaussian Process Classifier (GPC), a non-parametric model capable of capturing complex relationships in the
data~\cite[Sec.~6.4]{bishop_06}, that returns the estimated probability of Class~1—i.e., passing the time-domain specifications—so that $\hat f(\boldsymbol{\theta})\in[0,1]$.

The Class 1 probability at new input $\ve \theta$ given observed $\ve \theta_\mathcal{D}$ with labels $\ve \ell_\mathcal{D}$ is calculated by
\begin{align*}
p(\ell = 1 | \ve \theta, \boldsymbol{\theta}_\mathcal{D}, \boldsymbol{\ell}_\mathcal{D}) = \sigma\big(\mu_a \big(1+ \pi\sigma_a^2/8 \big)^{-1/2} \big),
\end{align*}
where $\sigma(a) \simeq \Phi(\lambda a)$ is the inverse probit function, while $\mu_a$ and $\sigma_a^2$ denote the posterior mean and variance of the latent variable at $\boldsymbol{\theta}$ under the Laplace approximation \cite[Sec.~3.4]{Rasmussen-Williams-06}, using data $\mathcal{D}_{\textnormal{hist}}$ from the historian together with newly collected experiments, $\mathcal{D}_{\textnormal{online}}$.

In short, $\hat f(\boldsymbol{\theta})$ provides the score we use to prioritize candidate controllers. The next step details how we select safe experiments via an acquisition function within a data-supported neighborhood around historical operation.

\subsection*{Smart sampling: Expected Improvement}
Given the previous observations and the constructed surrogate model, we adaptively pick new parameters within the trust region defined above to keep trials as minimal as possible. To this end, we use \emph{Expected Improvement (EI)}, a Bayesian optimization acquisition function that selects the most informative sampling point.

EI measures the expected reduction of the objective relative to the current best at a new point $\boldsymbol{\theta}\in\Theta$, i.e.,
\begin{align}
\text{EI}_n(\ve \theta) :&= \mathbb{E}_n\left\{[J_n^\ast - J(\ve \theta)]_+\right\} \notag\\
&= \mathbb{E}_n\big\{[d(\ve \theta^\ast, \ve \theta_0) + \lambda |\hat{f}(\ve \theta^\ast)-0.5|+\Omega\left(\ve \theta^*\right) \notag\\
&\hspace{2.15em} - d(\ve \theta, \ve \theta_0) - \lambda |\hat{f}(\ve \theta)-0.5|-\Omega(\ve \theta)]_+\big\}  \label{eq:EIn},
\end{align}
where $\mathbb{E}_n$ denotes expectation conditioned on the previously observed data, $\Omega(\cdot)$ represents the sparsity and plausibility terms in \eqref{eq:asce_cost}, $[a]_+ := \text{max}\{ a, 0\}$ and $\ve \theta^\ast = \text{argmin}_{\ve\theta_i \in \boldsymbol{\theta}_{1:n}} J(\ve\theta_i)$, is the cost function minimizer among the previously observed inputs.

To evaluate \eqref{eq:EIn}, ACE uses a Monte Carlo scheme that draws samples from the GPC-based model and approximates $\mathrm{EI}_n$ by averaging over those samples. Finally, we approximate the maximizer of $\mathrm{EI}_n$ with a general-purpose optimizer such as L-BFGS-B~\cite{nocedal_89}. Each EI proposal $\boldsymbol{\theta}_k$ is validated by a short closed-loop step test, after which we assign the pass/fail label $\ell_k=h(\boldsymbol{\theta}_k)\in\{0,1\}$. The online dataset is then updated, i.e., $\mathcal{D}_{\mathrm{online}}\leftarrow \mathcal{D}_{\mathrm{online}}\cup\{(\boldsymbol{\theta}_k,\ell_k)\}$, and the classifier $\hat f$ is retrained on the combined dataset $\mathcal{D}=\mathcal{D}_{\mathrm{hist}}\cup\mathcal{D}_{\mathrm{online}}$.

\begin{algorithm}
   \caption{ACE-Controller Fine Tuning}
   \label{alg:ace_dcs}

   \textbf{Input}: baseline controller parameters $\boldsymbol{\theta}_0$, historian archive $\mathcal{H}=\{(\boldsymbol{r}_1^{(i)},\boldsymbol{u}_1^{(i)},\boldsymbol{u}_2^{(i)},\boldsymbol{y}_1^{(i)},\boldsymbol{y}_2^{(i)})\}_{i=1}^{N_{\mathrm{hist}}}$ \\
   \textbf{Parameters}: $\lambda_0$ (init. penalty), $\lambda_{\max}$ (max penalty), $\epsilon$ (tolerance), $p$ (penalty growth) \\
   \textbf{Output}: CFE $\boldsymbol{\theta}_s^n$
   \begin{algorithmic}
     \STATE $\boldsymbol{\theta}^{(1:N_{\mathrm{hist}})} \leftarrow \textsc{IdentifyFromDCSHistorian}(\mathcal{H})$
  \STATE $\boldsymbol{\ell}_{1:N_{\mathrm{hist}}} \leftarrow \textsc{ClosedLoopMetrics}(\mathcal{H})$ \ 
  \STATE $\mathcal{D} \leftarrow \{(\boldsymbol{\theta}^{(i)},\ell_i)\}_{i=1}^{N_{\mathrm{hist}}}$;\ \ $\hat f \leftarrow \textsc{GPC\_Fit}(\mathcal{D})$

   \WHILE{$h(\boldsymbol{\theta}_0)=h(\boldsymbol{\theta}_s^n)$ \textbf{and} $\|\boldsymbol{\theta}_s^n-\boldsymbol{\theta}_0\|<\|\boldsymbol{\theta}_s^o-\boldsymbol{\theta}_0\|$}
     \STATE $\boldsymbol{\theta}_s^o \leftarrow \boldsymbol{\theta}_s^n$, \quad $\lambda_k \leftarrow \lambda_0$
     \WHILE{$\|\boldsymbol{\theta}_k-\boldsymbol{\theta}_{k-1}\|>\epsilon$ \textbf{or} $\lambda_k<\lambda_{\max}$}
       \STATE $\boldsymbol{\theta}_k \leftarrow \textsc{ExpectedImprovement}(\mathcal{D},\hat f, \ve \theta_0,\lambda_k)$
       \STATE $\ell_k \leftarrow \textsc{ClosedLoopStepTest}(\boldsymbol{\theta}_k)$ 
       \STATE $\mathcal{D} \leftarrow \mathcal{D}\cup\{(\boldsymbol{\theta}_k,{\ell}_k)\}$; \; $\hat f \leftarrow \textsc{GPC\_Update}(\mathcal{D})$
       \STATE $k \leftarrow k+1$, \quad $\lambda_k \leftarrow (\lambda_{k-1})^{p}$
     \ENDWHILE
     \STATE $\boldsymbol{\theta}_s^n \leftarrow \textsc{SampleDecisionBoundary}(\mathcal{D},\hat f,\boldsymbol{\theta}_0)$
   \ENDWHILE

   \STATE \textbf{return} $\boldsymbol{\theta}_s^n$
\end{algorithmic}
\end{algorithm}

\subsection{Algorithm}\label{sec:algo}

We present the full sample-efficient controller retuning algorithm based on Bayesian optimization and counterfactual explanations in Algorithm~\ref{alg:ace_dcs}. The method is initialized from the DCS archive: \textsc{IdentifyFromDCSHistorian} identifies controller parameters, \textsc{ClosedLoopMetrics} assigns pass/fail labels, and a GPC prior $\hat f$ is fitted to the resulting dataset $\{(\boldsymbol{\theta}^{(i)},\ell_i)\}_{i=1}^{N_{\mathrm{hist}}}$, where $\ell_i = h(\boldsymbol{\theta}^{(i)})$. Then \textsc{ExpectedImprovement} proposes actionable, trust-region candidates; each candidate is evaluated via a short closed-loop step test (\textsc{ClosedLoopStepTest}), after which the dataset and $\hat f$ are updated, and the boundary penalty parameter $\lambda_k$ in~\eqref{eq:asce_cost} is increased to steer queries toward $f(\boldsymbol{\theta})=0.5$. Finally, \textsc{SampleDecisionBoundary} selects a counterfactual near $f(\boldsymbol{\theta})=0.5$ with minimal proximity to $\boldsymbol{\theta}_0$, while sparsity and plausibility are enforced through the cost in~\eqref{eq:asce_cost}; here, $\boldsymbol{\theta}_s^n$ is the current CFE candidate and $\boldsymbol{\theta}_s^o$ the previously accepted one.

\section{Results and discussions}
\label{sec:results}

\begin{figure}
\centering
\includegraphics[trim=0mm 0mm 0mm 0mm,clip, width=0.95\columnwidth]{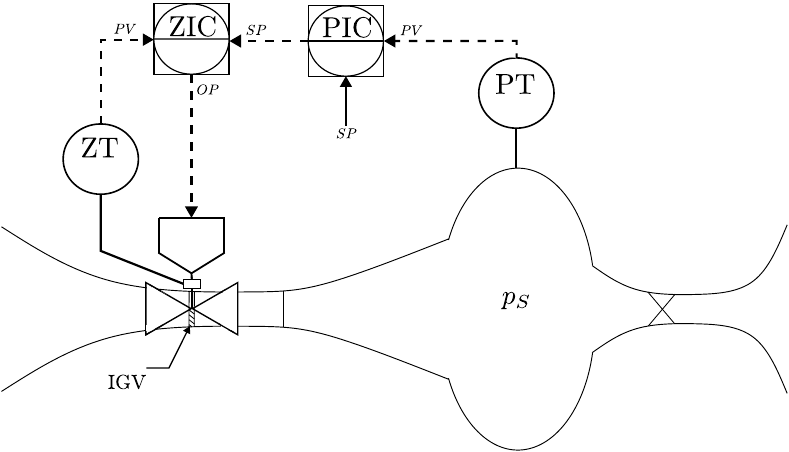}
\caption{P\&ID of the compressor system.}
    \label{fig:piid} 
\end{figure}

In this section, we evaluate Algorithm~\ref{alg:ace_dcs} on the compressor–plenum–throttle system. Figure~\ref{fig:piid} shows the site Piping and Instrumentation Diagram (P\&ID): the outer pressure controller (PIC) and the inner IGV position controller (ZIC) are the loops ACE retunes. We report results for two scenarios: \textbf{Case 1} (single-loop) in which IGV position feedback is unavailable and the PIC output drives the servo valve directly—hence ACE retunes only the outer loop; and \textbf{Case 2} (cascade), consistent with the P\&ID schematic, where IGV position feedback is present, the inner loop (ZIC) closes the vane-position dynamics, and ACE jointly retunes both controllers.


\subsection{Experimental setting}

\textbf{Controllers.} In both scenarios the controllers are implemented as PIDs with derivative filtering, i.e.,
\[
C(s)\;=\;K_p\;+\;\frac{K_i}{s}\;+\;\frac{K_d\, s}{1+T_f s}.
\]
Accordingly, ACE searches for a four-parameter retuning in \textbf{Case~1}, with $\boldsymbol{\theta}_{case_1}=[K_p, K_i, K_d, T_f]^\top$, and for an eight-parameter retuning in \textbf{Case~2}, with $\boldsymbol{\theta}_{case_2}=[K_{p1}, K_{i1}, K_{d1}, T_{f1}, K_{p2}, K_{i2}, K_{d2}, T_{f2}]^\top$, where the subscripts $1$ and $2$ on the PID parameters map to the outer ($C_1$) and inner ($C_2$) controllers, respectively; in both cases, ACE selects counterfactual retunings that satisfy the \emph{outer-loop} (PIC) closed-loop specifications.

\noindent \textbf{Parameters.} Unless otherwise noted, the compressor– plenum constants and IGV servo valve parameters are
$(\ell_c,B,k_T,a,\alpha)=(2.0,2.0,0.20,0.18,0.30)$ and $(K_v,\tau)=(1.0,0.5)$, respectively. We use a normalized step of amplitude $r_{\text{step}}=1$. Measurement noises are zero-mean Gaussian with variances
$(\sigma_{\text{outer}}^2,\sigma_{\text{inner}}^2)=(0.01,0.005)$, corresponding to a outer-loop input signal-to-noise ratio of $20[\textnormal{dB}]$. ACE hyperparameters follow  \cite[Appendix G.1]{Guerrero-25}.

\noindent \textbf{Continuous-to-Discrete conversion.} The compressor–plenum path and the servo valve are mapped from $G_1(s)$ and $G_2(s)$ to $G_1^{\textnormal{d}}(q)$ and $G_2^{\textnormal{d}}(q)$ using the bilinear (Tustin) transform with sampling period $T_s=0.1\,\text{s}$. The PID controllers are implemented in position form: the derivative term uses a first–order filter discretized with Tustin–consistent coefficients, while the integral term uses the rectangular (forward Euler) rule.

\noindent \textbf{Closed-loop criteria.}
We define the pass/fail function as
\[
h(\boldsymbol{\theta})=
\begin{cases}
1, &
\text{if }\hspace{0.3em}
\begin{aligned}
& \hspace{1.5em} |e_{\mathrm{ss}}| \le 0.01 \\
& \&\;\hspace{0.5em} t_{\mathrm{rise}\,(10\to90)}\le 20~\text{s} \\
& \&\;\hspace{0.5em} t_{\mathrm{settle}} \le 50~\text{s} \\
& \&\;\hspace{0.5em} \mathrm{OS} \le 20\%
\end{aligned}
\\[1ex]
0, & \hspace{1em}\text{otherwise,}
\end{cases}
\]
where $e_{\mathrm{ss}}$ is the steady-state error, $t_{\mathrm{rise}(10\to90)}$ is the rise time between the first crossings of $10\%$ and $90\%$ of the final step change, $t_{\mathrm{settle}}$ is the first time for which $|y(t)-y_\infty|\le 0.05$ for all $t\ge t_{\mathrm{settle}}$, and $\mathrm{OS}$ is the percent of overshoot.



\subsection{Results}
\subsection*{ Case 1}

Table~\ref{tab:case1} summarizes results over 100 ACE runs for three fixed instances $\boldsymbol{\theta}_{01:03}$, initialized with 30 historical controllers identified from the DCS archive. Each entry reports the mean \emph{Validity} (fraction of runs that found a counterfactual), the mean \emph{LOF} (Local
Outlier Factor; values close to or greater than $1$ indicate inliers to the data manifold) computed over every new
$\ve{\theta}_k$ proposed during the search (including the final CFE), the mean number of additional closed-loop tests (\emph{iters.}), and the mean component-wise parameter shift reported as $(K_p;K_i;K_d;T_f)$, where we define the shift as $\Delta=\boldsymbol{\theta}_{\text{CFE}}-\boldsymbol{\theta}_0$---i.e., positive means increase w.r.t. baseline. Overall, average changes are small (typically $|\Delta|\!\lesssim\!0.1$), with only a modest number of extra step tests per run. Over 100 trials, ACE requires $21.6$~s per valid CFE on average, demonstrating practical deployability.

\begin{table}[t]
    \centering
    \caption{Case 1 (single loop): summary over 100 ACE runs per instance.}
    \label{tab:case1}
    \scriptsize
    \setlength{\tabcolsep}{3pt}
    \resizebox{\columnwidth}{!}{
    \begin{tabular}{l|c|c|c|c}
        \hline
        Inst. & Val. & LOF & iters. &
         $\ve\Delta K_p; \ve\Delta K_i; \ve\Delta K_d; \ve\Delta T_f$ \\
        \hline
        $\boldsymbol{\theta}_{01}$ & 1 & 1.0276 & 12 &
        $+0.0334;\ -0.0986;\ +0.0251;\ +0.0331$ \\
        $\boldsymbol{\theta}_{02}$ & 1 & 1.0251 & 12  &
        $+0.0057;\ -0.0482;\ +0.0029;\ +0.0333$ \\
        $\boldsymbol{\theta}_{03}$ & 1 & 1.0205 & 18  &
        $+0.0234;\ -0.0602;\ +0.0071;\ +0.0325$ \\
        \hline
    \end{tabular}}
\end{table}

\begin{figure}
\centering
\includegraphics[trim=0mm 0mm 0mm 0mm,clip, width=1\columnwidth]{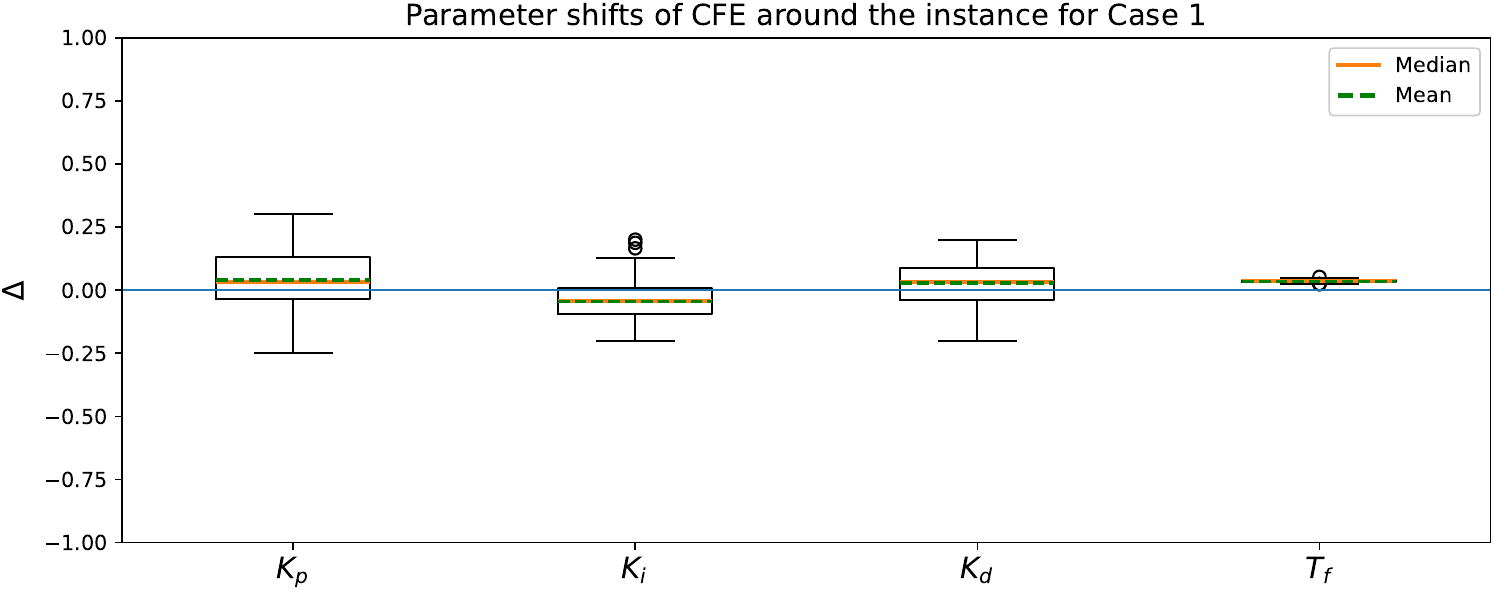}
\caption{Case 1 - Box plot of component-wise shifts $\ve \Delta=\boldsymbol{\theta}_{\text{cand}}-\boldsymbol{\theta}_0$ across 100 ACE runs.}
    \label{fig:box_plot_single} 
\end{figure}
\begin{figure}
\centering
\includegraphics[trim=0mm 0mm 0mm 0mm,clip, width=1\columnwidth]{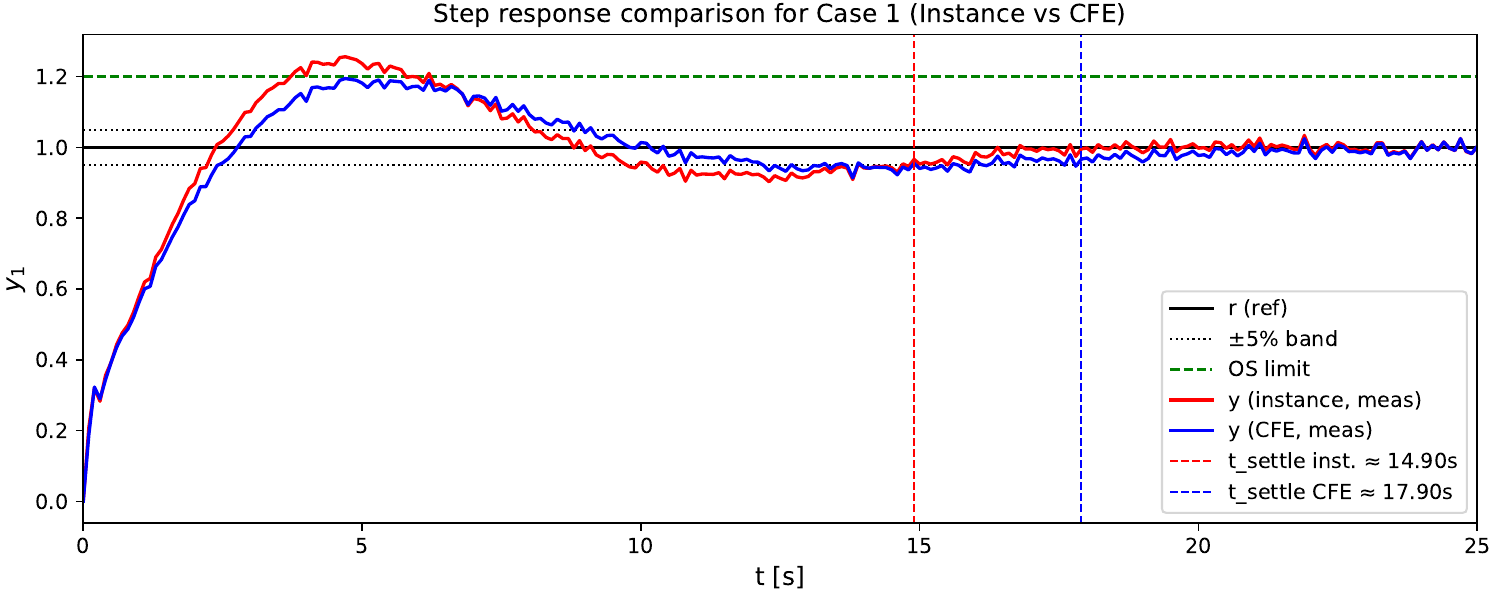}
\caption{Case 1 - Step responses: baseline $\boldsymbol{\theta}_0$ (red) vs. one representative CFE (blue).}
    \label{fig:step_one} 
\end{figure}
To examine the results in detail, we focus on the first instance in Table~\ref{tab:case1}: $\boldsymbol{\theta}_0= \ve\theta_{01}=[3.864,1.112,1.561,0.015]$. Figure~\ref{fig:box_plot_single} reports the component-wise shifts $\Delta=\boldsymbol{\theta}_{\text{cand}}-\boldsymbol{\theta}_0$ across 100 ACE runs—including, but not limited to, the final CFEs. Note that, unlike Table~\ref{tab:case1}—which averages only the \emph{final} CFE per run—Figure~\ref{fig:box_plot_single} pools \emph{all} candidate evaluations generated during the ACE search, with an average of $\sim$12 proposals per run until the CFE is found. The boxes span the interquartile range (25–75\%), the solid orange line indicates the median and the dashed green line indicates the mean,  whiskers extend to the non-outlier range, and circles denote outliers. Three patterns emerge. First, the distributions are tightly centered around zero, with few outliers and typical magnitudes $|\Delta|\lesssim 0.3$, indicating that ACE proposes small, local retunings consistent with the trust-region constraint. Second, the medians suggest a mild \emph{increase} in $K_p$, a mild \emph{decrease} in $K_i$, and small positive shifts for $K_d$ and $T_f$—in line with making the loop slightly more responsive while filtering derivative action. Third, several CFEs exhibit $\Delta \approx 0$ on individual coordinates, showing that sparsity is naturally promoted: not every retuning needs to change all gains.

Figure~\ref{fig:step_one} compares the instance $\boldsymbol{\theta}_0=[3.864,1.112,1.561,$ $0.015]$ (red) against one representative counterfactual $\mathrm{CFE}=[3.916,1.007,1.542,0.022]$ (blue). Consistent with Figure~\ref{fig:box_plot_single}, the counterfactual increases $K_p$ slightly, reduces $K_i$, and applies modest adjustments to $K_d$ and $T_f$. These small parameter changes are sufficient to flip the label from “fail” to “pass”: overshoot is brought below the 20\% limit (peak $\le 1.20$), while all other time–domain metrics are satisfied. The comparison illustrates the core CFE objective: achieving specification compliance with minimal movement in parameter space and without leaving the safe operating neighborhood; from a control perspective, this is achieved through small local gain adjustments that reduce overshoot and improve transient compliance without altering the closed-loop operating regime.

\subsection*{ Case 2}
As in Case~1, Table~\ref{tab:case2} summarizes 100 ACE runs for three fixed cascade instances—\emph{Validity} is $1.0$ for all instances and the mean \emph{LOF} values remain close to $1$, indicating that the proposed controllers lie on (or near) the historical manifold. Despite doubling the decision variables (eight gains instead of four), the mean number of additional closed-loop tests per run remains modest—$\sim$15–18, comparable to Case~1—showing that ACE can jointly retune inner and outer loops without the usual two-step manual procedure used in industry~\cite[Chap.~16, p.~284]{seborg2016}. Parameter shifts remain small in magnitude (typically $|\Delta|\!\lesssim\!0.4$). Over 100 trials, ACE requires $29.5$~s per valid CFE on average. In a sensitivity check, varying the historian archive size $\mathcal{H}$ yields broadly stable metrics in both cases, indicating a robust trust-region initialization. No closed-loop instability was observed for any candidate $\ve{\theta}_k$ evaluated in this study.
\begin{table}
    \centering
    \caption{Case 2 (cascade): summary over 100 ACE runs per instance.}
    \label{tab:case2}
    \scriptsize
    \setlength{\tabcolsep}{3pt}
    {\renewcommand{\arraystretch}{0.9}
    \resizebox{\columnwidth}{!}{
    \begin{tabular}{l|c|c|c|c}
        \hline
        Inst. & Val. & LOF & iters. &
        \shortstack[c]{\textbf{$C_1,C_2$:} $\ve\Delta K_p; \ve\Delta K_i; \ve\Delta K_d; \ve\Delta T_f$} \\
        \hline
        \raisebox{0.6ex}{$\boldsymbol{\theta}_{01}$} & \raisebox{0.6ex}{1} & \raisebox{0.6ex}{1.0218} & \raisebox{0.6ex}{15} &
        \shortstack[c]{%
            \rule{0pt}{2ex}$+0.0488;\ +0.3120;\ +0.0276;\ +0.0418$\\[-0.9ex]%
            $+0.0226;\ -0.0839;\ +0.0587;\ +0.0302$%
        } \\
        \raisebox{0.6ex}{$\boldsymbol{\theta}_{02}$} & \raisebox{0.6ex}{1} & \raisebox{0.6ex}{1.0123} & \raisebox{0.6ex}{17} &
        \shortstack[c]{%
            \rule{0pt}{2ex}$+0.1106;\ +0.1144;\ +0.2259;\ +0.0404$\\[-0.9ex]%
            $-0.1123;\ -0.0432;\ +0.1088;\ +0.0415$%
        } \\
        \raisebox{0.6ex}{$\boldsymbol{\theta}_{03}$} & \raisebox{0.6ex}{1} & \raisebox{0.6ex}{1.0285} & \raisebox{0.6ex}{18} &
        \shortstack[c]{%
            \rule{0pt}{2ex}$+0.3958;\ +0.1374;\ +0.0780;\ +0.0391$\\[-0.9ex]%
            $+0.0546;\ -0.1348;\ +0.0614;\ +0.0405$%
        } \\
        \hline
    \end{tabular}}}
\end{table}

\begin{figure}
\centering
\includegraphics[trim=0mm 0mm 0mm 0mm,clip, width=1\columnwidth]{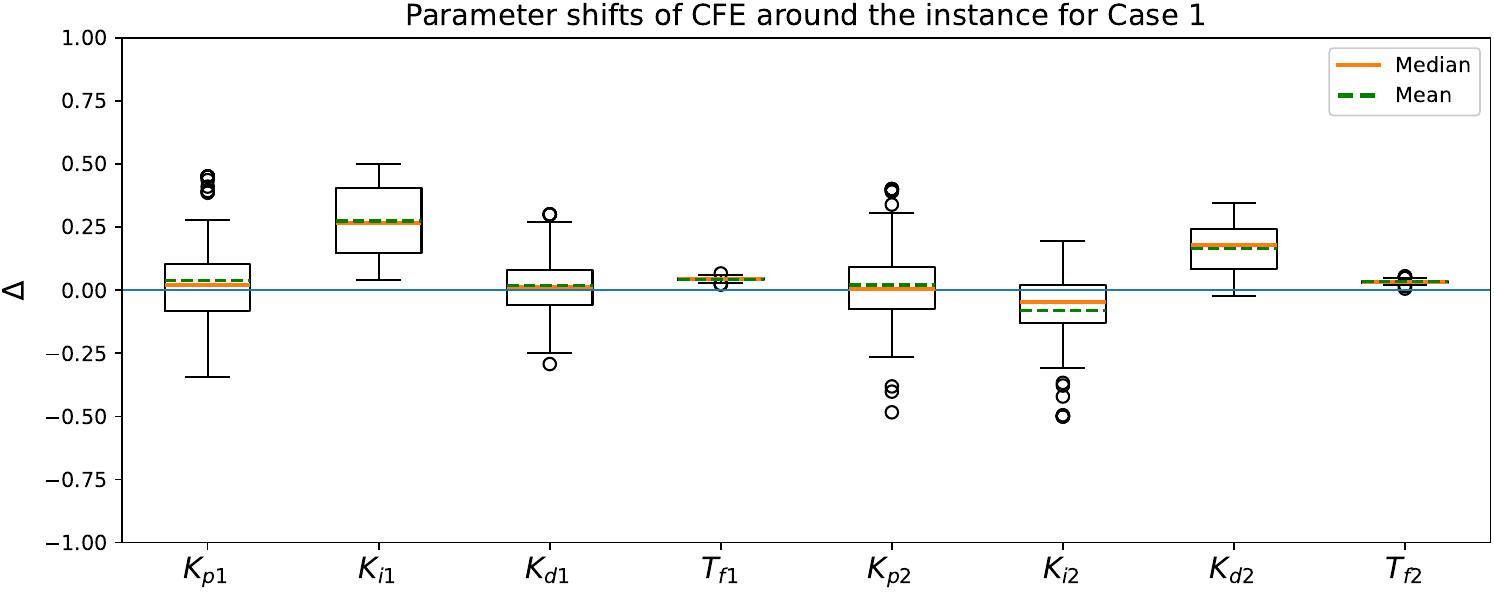}
\caption{Case 2 - Box plot of component-wise shifts $\ve \Delta=\boldsymbol{\theta}_{\text{cand}}-\boldsymbol{\theta}_0$ across 100 ACE runs.}
    \label{fig:box_plot_dual} 
\end{figure}

Figure~\ref{fig:box_plot_dual} summarizes the component-wise shifts
$\Delta=\boldsymbol{\theta}_{\text{cand}}-\boldsymbol{\theta}_0$
across 100 ACE runs for instance $\boldsymbol{\theta}_{01}$. As in Case~1, the distributions remain tightly
centered around zero with few outliers, indicating small, local moves
consistent with the trust region. Two coordinates show a clear directional
trend: $K_{i1}$ and $K_{d2}$ tend to increase, both typically within $|\Delta|\lesssim 0.5$. In contrast,
$K_{i2}$ exhibits a slight decrease, while $T_{f1}$ and $T_{f2}$ show mild
positive shifts; the remaining gains are essentially centered at
$\Delta\approx 0$. Overall, even in the 8-parameter search, ACE stays near
the historical manifold and proposes sparse, low-magnitude adjustments.

\begin{figure}
\centering
\includegraphics[trim=0mm 0mm 0mm 0mm,clip, width=1\columnwidth]{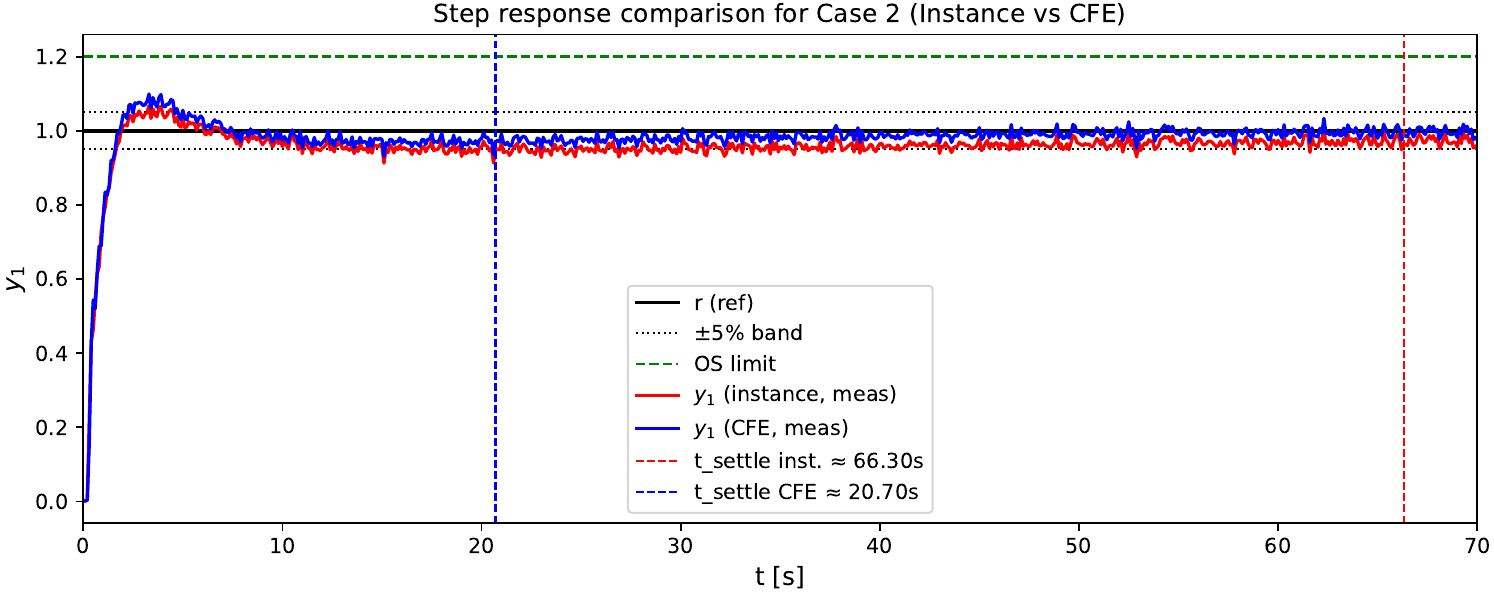}
\caption{Case 2 - Step responses: baseline $\boldsymbol{\theta}_0$ (red) vs. one representative CFE (blue).}
    \label{fig:step_dual} 
\end{figure}

A representative comparison is shown in Fig.~\ref{fig:step_dual}: the
instance $\boldsymbol{\theta}_0=[8.453,0.136,1.73,0.010,3.986,2.842,0.264$, $0.02]$
versus one CFE $[8.519,0.466,1.73,$$0.065,3.84,2.64,$ $0.569,0.065]$. The parameter change is $\Delta_{\text{CFE}}=\text{CFE}-\boldsymbol{\theta}_0 =[+0.066, +0.331, +0.003, +0.055,$ $-0.146,-0.202$, $+0.306,+0.045]$,
i.e., a modest boost in $K_{i1}$ and $K_{d2}$, slight increases in
$T_{f1}$/$T_{f2}$, and small shifts elsewhere. In PID terms, the larger
$K_{i1}$ removes steady-state error more decisively in the outer loop, while
the higher $K_{d2}$ improves inner-loop vane tracking and effective damping.
The result is a substantial reduction of settling time—from $\approx 66.3$\,s
down to $\approx 20.7$\,s—while keeping the overshoot within limits and
satisfying all four time-domain metrics, achieved with small, interpretable
moves in parameter space.

\section{Conclusions}
\label{sec:conclusions}
We have presented a historian–anchored, counterfactual explainability approach to safe, sample-efficient controller retuning, demonstrated on the axial-compressor pressure loop of a compressor–plenum–throttle system. Our ACE-style method is model-free, operates entirely in closed loop, and proposes minimal, actionable parameter changes within a trust region, requiring only a handful of short closed-loop step tests. In both single-loop and cascade architectures, ACE produced small, interpretable shifts in PID gains that flipped the pass/fail label from “fail” to “pass” while staying on the historical data manifold and preserving surge margin. Overall, ACE delivers sample-efficient, low-impact retuning while maintaining safety margins by remaining within
the historian-supported operating envelope. Beyond compressors, the method can be envisioned for data-supported, costly-query retuning, and may be conditioned on operating context (including time) to address nonstationary regimes. Future work includes integrating instrumental variable and kernel-regularized estimators to refine historian-based priors for general closed-loop configurations, extending beyond PID to nonlinear controllers, whose more complex control laws may overcome fundamental limitations of linear control but entail more sophisticated performance metrics and robustness considerations, and validating the approach on real DCS historian data and/or high-fidelity industrial simulators.
\section*{Acknowledgements}
The authors thank Ricardo A. Guerrero Alvarez for his valuable insight into compressor systems and for sharing practical knowledge that contributed to this work.

\bibliography{references} 

@book{Rasmussen-Williams-06,
  author    = {C. E. Rasmussen and C. K. I. Williams},
  publisher = {MIT Press},
  title     = {Gaussian Processes for Machine Learning},
  year      = {2006},
}

@techreport{moore1985theory,
  title={A theory of post-stall transients in multistage axial compression systems},
  author={Moore, F. K. and Greitzer, E. M.},
  year={1985},
  institution={NASA}
}

@book{Ljung1999SystemUser,
    title = {{System Identification: Theory For The User}},
    year = {1999},
    author = {Ljung, L.},
    edition = {Second},
    publisher = {Prentice-Hall}
}

@book{Soderstrom2001SystemIdentification,
    title = {{System Identification}},
    year = {2001},
    author = {S{\"{o}}derstr{\"{o}}m, T. and Stoica, P.},
    publisher = {Prentice Hall}
}

@book{aastrom1984computer,
	title={Computer {C}ontrolled {S}ystems: {T}heory and {D}esign},
	author={{\AA}str{\"o}m, K. J. and Wittenmark, B.},
	publisher={Prentice-Hall},
	year={1984}
}

@article{hjalmarsson1998iterative,
  title={Iterative feedback tuning: theory and applications},
  author={Hjalmarsson, H. and Gevers, M. and Gunnarsson, S. and Lequin, O.},
  journal={IEEE Control Systems Magazine},
  volume={18},
  number={4},
  pages={26--41},
  year={1998}
}

@inproceedings{anderson1991adaptive,
  title={{Adaptive robust control: Online learning}},
  author={Anderson, B. D. O. and Kosut, R. L.},
  booktitle={30th IEEE Conference on Decision and Control},
  pages={297--298},
  year={1991}
}

@book{albertos2012iterative,
  title={Iterative identification and control: Advances in theory and applications},
  author={Albertos, P. and Piqueras, A. S.},
  year={2012},
  publisher={Springer}
}

@article{van1995identification,
  title={Identification and control—{Closed}-loop issues},
  author={Van Den Hof, P. M. J. and Schrama, R. J. P.},
  journal={Automatica},
  volume={31},
  number={12},
  pages={1751--1770},
  year={1995}
}

@article{campi2002virtual,
  title={Virtual reference feedback tuning: a direct method for the design of feedback controllers},
  author={Campi, M. C. and Lecchini, A. and Savaresi, S. M.},
  journal={Automatica},
  volume={38},
  number={8},
  pages={1337--1346},
  year={2002},
  publisher={Elsevier}
}

@book{aastrom2006advanced,
  title={Advanced PID Control},
  author={{\AA}str{\"o}m, K. J. and H{\"a}gglund, T.},
  year={2006},
  publisher={ISA-The Instrumentation, Systems and Automation Society}
}

@book{Molnar_22,
  author  = {Molnar, C.},
  title   = {Interpretable machine learning: A guide for making Black Box models explainable},
  year    = 2022,
  publisher = {Independently published}
}

@article{wachter_17,
  author  = "Wachter, Sandra and Mittelstadt, Brent and Russell, Chris",
  title   = {Counterfactual Explanations Without Opening the Black Box: Automated Decisions and the {GDPR}},
  journal = "Harvard Journal of Law \& Technology",
  year    = 2017,
  volume  = 31,
  pages   = "841-887",
}

@inproceedings{allamaa2025,
  title={{ExAMPC: The data-driven explainable and approximate NMPC with physical insights}},
  author={Allamaa, J. P. and Patrinos, P. and Son, T. D.},
  booktitle={2025 IEEE/RSJ International Conference on Intelligent Robots and Systems (IROS)},
  pages={2976--2983},
  year={2025},
  organization={IEEE}
}

@inproceedings{porcari2025,
  title={{eXplainable AI for data driven control: an inverse optimal control approach}},
  author={Porcari, F. and Materassi, D. and Formentin, S.},
  booktitle={2025 IEEE 64th Conference on Decision and Control (CDC)},
  pages={5186--5191},
  year={2025},
}

@inproceedings{Ribeiro16,
  title={{Why should I trust you?" Explaining the predictions of any classifier}},
  author={Ribeiro, Marco Tulio and Singh, Sameer and Guestrin, Carlos},
  booktitle={{22nd ACM SIGKDD}},
  pages={1135--1144},
  year={2016}
}

@inproceedings{Lundberg17,
  title={A unified approach to interpreting model predictions},
  author={Lundberg, Scott M and Lee, Su-In},
  booktitle={31st International Conference on Neural Information Processing Systems},
  pages={4768--4777},
  year={2017}
}

@ARTICLE{Riva24,
  title={{Toward eXplainable data-driven control (XDDC): The property-preserving framework}},
  author={Riva, G. and Formentin, S.},
  journal={IEEE Control Systems Letters},
  volume={8},
  pages={478--483},
  year={2024}
}

@article{kurz2012gas,
  title={Gas compressor station economic optimization},
  author={Kurz, Rainer and Lubomirsky, Matt and Brun, Klaus},
  journal={International Journal of Rotating Machinery},
  volume={2012},
  number={1},
  pages={715017},
  year={2012},
  publisher={Wiley Online Library}
}

@book{bishop_06,
  author    = {C. M. Bishop},
  title     = {Pattern Recognition and Machine Learning},
  publisher = {Springer},
  year      = {2006},

}

@book{seborg2016,
  title={Process Dynamics and Control},
  author={Seborg, Dale E and Edgar, Thomas F and Mellichamp, Duncan A and Doyle III, Francis J},
  year={2016},
  publisher={John Wiley \& Sons}
}

@article{beagle21,
  title={Heavy-Duty Gas Turbine Operating and Maintenance Considerations},
  author={Beagle, D. and Moran, B. and McDufford, M. and Merine, M.},
  journal={GE Vernova Atlanta, GA},
  year={2021}
}

@article{smith2001,
  title={A review of air separation technologies and their integration with energy conversion processes},
  author={Smith, Arthur R and Klosek, Joseph},
  journal={Fuel Processing Technology},
  volume={70},
  number={2},
  pages={115--134},
  year={2001},
  publisher={Elsevier}
}

@incollection{crompton2021,
  title={Data Management from the {DCS} to the {Historian}},
  author={Crompton, J.},
  booktitle={Machine Learning and Data Science in the Oil and Gas Industry},
  editor={Patrick Bangert},
  pages={83--110},
  year={2021}
}

@inproceedings{jager95,
  title={Rotating stall and surge control: A survey},
  author={De Jager, Bram},
  booktitle={34th IEEE Conference on Decision and Control},
  volume={2},
  pages={1857--1862},
  year={1995},
}

@article{greitzer1976surge,
  title={Surge and rotating stall in axial flow compressors: Theoretical compression system model},
  author={Greitzer, Edward M},
  journal={Journal of Engineering for Power},
  volume={98},
  number={2},
  pages={190--198},
  year={1976},
  publisher={American Society of Mechanical Engineers Digital Collection}
}

@article{brandt1994,
  title={{GE} gas turbine design philosophy},
  author={Brandt, DE and Wesorick, RR},
  journal={GER-3434, General Electric},
  year={1994}
}

@inproceedings{breunig_00,
author = {Breunig, M. M. and Kriegel, H.-P. and Ng, R. T. and Sander, J.},
title = {{LOF}: identifying density-based local outliers},
year = {2000},
booktitle = {{2000 ACM SIGMOD}},
pages = {93–104},
}

@inproceedings{Guerrero-25,
title = {{ACE}: Adapting Sampling for Counterfactual Explanations},
author = {Guerrero, M. A. and Rojas, C. R.},
booktitle = {8th Annual Learning for Dynamics \& Control Conference},
year = {2026},
publisher= {Also available at arXiv:2509.26322},
}

@article{nocedal_89,
  author    = {Dong C. Liu and Jorge Nocedal},
  title     = {On the limited memory {BFGS} method for large scale optimization},
  journal   = {Mathematical Programming},
  year      = {1989},
  volume    = {45},
  pages     = {503--528},
}
\end{document}